\documentclass[10pt]{article}
\usepackage{amsmath,amssymb,amsfonts}

\usepackage{graphicx} 

\newtheorem{theorem}{Theorem}
\newtheorem{lemma}{Lemma}
\newtheorem{definition}{Definition}

\date{}

\title{The Limiting Distribution of the Hook Length of a Randomly
Chosen Cell in a Random Young Diagram}
\bigskip

\author{{\bf Ljuben Mutafchiev}\\
American University in Bulgaria, 2700 Blagoevgrad, Bulgaria \\ and
Institute of Mathematics and Informatics of the \\ Bulgarian
Academy of Sciences
\\ \tt {ljuben@aubg.bg}}
\begin{document} \maketitle

\begin{abstract}
Let $p(n)$ be the number of all integer partitions of the positive
integer $n$ and let $\lambda$ be a partition, selected uniformly
at random from among all such $p(n)$ partitions. It is known that
each partition $\lambda$ has a unique graphical representation,
composed by $n$ non-overlapping cells in the plane called Young
diagram. As a second step of our sampling experiment, we select a
cell $c$ uniformly at random from among the $n$ cells of the Young
diagram of the partition $\lambda$. For large $n$, we study the
asymptotic behavior of the hook length $Z_n=Z_n(\lambda,c)$ of the
cell $c$ of a random partition $\lambda$. This two-step sampling
procedure suggests a product probability measure, which assigns
the probability $1/np(n)$ to each pair $(\lambda,c)$. With respect
to this probability measure, we show that the random variable $\pi
Z_n/\sqrt{6n}$ converges weakly, as $n\to\infty$, to a random
variable whose probability density function equals
$6y/\pi^2(e^y-1)$ if $0<y<\infty$, and zero elsewhere.
\end{abstract}

\vspace{.5cm}



\vspace{.2cm}

\section{Introduction and Statement of the Main Result}

For a natural number $n$, we say that $\lambda$ is a partition of
$n$ if $\lambda$ is a sequence
$(\lambda_1,\lambda_2,...,\lambda_k)$ of positive integers
satisfying $\lambda_1\ge\lambda_2\ge...\ge\lambda_k$ and such that
\begin{equation}\label{sigma}
\sum_{j=1}^k\lambda_j=n.
\end{equation}
The summands $\lambda_j$ in (\ref{sigma}) are usually called parts
of $\lambda$. The Young diagram of a partition is an array of
square boxes, or cells, in the first quadrant of the plane,
left-justified, with $\lambda_t$ cells in the $t$-th row counting
from the bottom. We label these cells $(t,s)$, with $t$ denoting
the row number of the cell and $s$ - the column number in the
Young diagram. For example, in the Young diagram of
$\tilde{\lambda}=(5,4,3,3,2,2,2,1)$ of the partition of $n=22$ as
$22=5+4+3+3+2+2+2+1$, the cell $(3,2)$ is the square whose
vertices in the $(s,t)$-plane have coordinates $(1,2), (2,2),
(2,3), (1,3)$. Reading consecutively the numbers of cells in the
columns of the array of the partition $\lambda$ , beginning from
the most left column, we get the conjugate partition
$\lambda^*=(\lambda_1^*,\lambda_2^*,...,\lambda_l^*)$, where
$l=\lambda_1$. The Young diagram of $\lambda^*$ is called
conjugate of the Young diagram of $\lambda$. In our example, we
have $\tilde{\lambda}^*=(8,7,4,2,1)$, which is the conjugate
partition of $22$, namely, $22=8+7+4+2+1$. The hook length of the
cell $c=(t,s)$ in the partition $\lambda$ is defined by the
formula
\begin{equation}\label{hook}
h(\lambda,c):=\lambda_t-s+\lambda_s^*-t+1,
\end{equation}
that is, $h(\lambda,c)$ is the number of cells in the hook
comprised by the $(t,s)$-cell itself and by the cells in the
$t$-th row right of $(t,s)$ and $s$-th column above $(t,s)$.
Recalling our example given above, we get
$h(\tilde{\lambda},(3,2))=6$ since $t=3,\lambda_3=3,s=2$ and
$\lambda_2^*=7$.

Let $\Lambda(n)$ be the set of all partitions of $n$ and let
$p(n)=|\Lambda(n)|$ (further on, by $|T|$  we denote the
cardinality of the set $T$). The number $p(n)$ is determined
asymptotically by the famous partition formula of Hardy and
Ramanujan \cite{HR18}:
\begin{equation}\label{hr}
p(n)\sim\frac{1}{4n\sqrt{3}}
\exp{\left(\pi\sqrt{\frac{2n}{3}}\right)}, \quad n\to\infty.
\end{equation}
A precise asymptotic expansion for $p(n)$ was found later by
Rademacher \cite{R37} (more details may be also found in
\cite[Chapter 5]{A76}). Further on, we assume that, for fixed
integer $n\ge 1$, a partition $\lambda\in\Lambda(n)$ is selected
uniformly at random (uar). In other words, we assign the
probability $1/p(n)$ to each $\lambda\in\Lambda(n)$. In this way,
each numerical characteristic of $\lambda$ can be regarded as a
random variable defined on $\Lambda(n)$ (or, a statistic in the
sense of the random generation of partitions of $n$). The study of
the asymptotic behavior of various partition statistics for large
$n$ is a subject of intensive research in combinatorics, number
theory and statistical physics. Erd\"{o}s and Lehner \cite{EL41}
were apparently the first who established a probabilistic limit
theorem related to integer partitions. As a matter of fact, they
found an appropriate normalization for the number of parts in a
random partition of $n$ and showed that it converges weakly to the
extreme value (Gumbel) distribution as $n\to\infty$. By
conjunction of the corresponding Young diagram, the same limit
theorem holds true for the largest part of a random integer
partition. For other typical distributional results and limit
theorems of various integer partition statistics, we refer the
reader, e.g., to \cite{ST771}, \cite{ST772}, \cite{ST78},
\cite{ES84}, \cite{F93}, \cite{P97}, \cite{M05}, \cite{W11},
\cite{GKW14}.

Another subject of study in the asymptotic theory of integer
partitions is related to the limit shape of the underlying Young
diagrams. Here is a simple setting of the problem. Let
\begin{equation}\label{multi}
l_j=l_j(\lambda):=|\{q:\lambda_q=j\}|
\end{equation}
 be the multiplicity of the
part equal to $j$ in a partition $\lambda\in\Lambda(n)$,
$j=1,2,...,n$. The upper boundary of the Young diagram of
$\lambda$ is a piecewise constant function
$X_\lambda:[0,\infty)\to\mathbb{Z}_+:=\{0,1,...\}$ given by
$$
X_\lambda(t):=\sum_{j\ge t} l_j.
$$
If the set $\Lambda(n)$ is endowed with a probability measure
$\mu_n$ (e.g., $\mu_n$ is the uniform measure $\mathbb{P}$ or the
Plancherel distribution on the set of Young diagrams with $n$
cells), the limit shape, with respect to $\mu_n$ as $n\to\infty$,
is understood as a function (curve) $s=s(t)$ in the plane $(t,s)$
such that for every $\delta>0$ and any $\epsilon>0$,
\begin{equation}\label{limit}
\lim_{n\to\infty}\mu_n \{\lambda\in\Lambda(n):
\sup_{t\ge\delta}{|A_n^{-1}X_\lambda(tB_n)-s(t)|}>\epsilon\}=0,
\end{equation}
where $A_n$ and $B_n$ are suitable scaling constants satisfying
$A_n B_n=1$. With respect to the uniform distribution $\mathbb{P}$
on $\Lambda(n)$ the limit shape in (\ref{limit}) exists under the
scaling $A_n=B_n=\sqrt{n}$ and is determined by the function
\begin{equation}\label{shapefunc}
s(t)=-\log{(1-e^{-\pi t/\sqrt{6}})},
\end{equation}
or, in a more symmetric form, by the equation
\begin{equation}\label{shapecur}
e^{-\pi s/\sqrt{6}}+e^{-\pi t/\sqrt{6}}=1.
\end{equation}
The limit shape (\ref{shapecur}) for the uniform distribution on
$\Lambda(n)$ was first identified by Temperley \cite{T52} in
relation to a model of a growing crystal. Afterward Szalay and
Tur\'{a}n \cite{ST771} obtained essentially analogous estimates of
those in (\ref{limit}) and (\ref{shapefunc}), which were used
later by Vershik \cite{V96} to establish (\ref{shapecur}) and
generalize it to other types of probability measures and models.
In the same paper \cite{V96}, Vershik gave also an important
probabilistic frame of the relationship between problems from
statistical mechanics (models of ideal gas) and problems related
to the asymptotic theory of partitions. An alternative proof of
(\ref{shapecur}) was also given by Pittel in \cite{P97}. Recent
developments in this area are presented in \cite{B14}.

 In the present paper we study a statistic
 produced by a random selection of a cell in a random Young diagram.
 The sampling procedure combines the outcomes of two experiments:
 first, we select a partition $\lambda\in\Lambda(n)$ uar, and
then we select a cell $c=(t,s)\in\lambda$ uar. A similar type of
sampling was first proposed and studied in the context of part
multiplicities of a random integer partition in \cite{CPSW99}, and
then in the context of part sizes in \cite{M151} and \cite{M152}.
The sampling procedure that we introduced leads to the product
probability space
\begin{equation}\label{omega}
\Omega(n)=\Lambda(n)\times\{c\in\lambda:\lambda\in\Lambda(n)\}
\end{equation}
 and
the product probability measure $\mathbb{P}(.)$, defined by
\begin{equation}\label{pm}
\mathbb{P}((\lambda,c)\in\Omega(n))=\frac{1}{np(n)}.
\end{equation}
(For more details on the theory of product probability spaces, we
refer the reader, e.g., to \cite[Part I, Section 4.2]{L60}). We
also denote by $\mathbb{E}(.)$ the expectation with respect to the
probability measure $\mathbb{P}$.

Our aim in this paper is to determine asymptotically, as
$n\to\infty$, the distribution of the hook length of a randomly
chosen cell of the Young diagram of a partition
$\lambda\in\Lambda(n)$ chosen uar. More precisely, for any pair
$(\lambda,c)\in\Omega(n)$, we denote the underlying hook length by
\begin{equation}\label{z}
 Z_n=Z_n(\lambda,c):=h(\lambda,c),
\end{equation}
 where the function $h(\lambda,c)$ was defined earlier by
(\ref{hook}). We organize the paper as follows. Section 2 contains
some auxiliary facts that we need further. In Section 3 we prove
the following limit theorem for $Z_n$.

\begin{theorem}
Let $0<y<\infty$. Then, we have
$$
lim_{n\to\infty}\mathbb{P}(\pi Z_n/\sqrt{6n}\le y)
=\frac{6}{\pi^2}\int_0^y\frac{u}{e^u-1} du.
$$
\end{theorem}

{\it Remark 1.} The hook lengths of cells in the Young diagram of
an integer partition play an important role in algebraic
combinatorics thanks to the famous hook-length formula due to
Frame et al. \cite{FRT54}. It represents the number of standard
Young tableaux $d(\lambda)$ of shape $\lambda\in\Lambda(n)$ by the
hook lengths of the cells $c\in\lambda$ in the following way:
\begin{equation}\label{tableaux}
d(\lambda)=\frac{n!}{\prod_{c\in\lambda}h(\lambda,c)}.
\end{equation}
A standard Young tableau of shape $\lambda$ is a labelling of the
cells of the Young diagram of $\lambda$ with numbers $1$ to $n$ so
that the labels are strictly increasing from bottom to top along
columns and from left to right along rows. There exists a
bijection between the partitions from the set $\Lambda(n)$ and the
irreducible representations of the symmetric group $S_n$ of $n$
letters (see, e.g., \cite[Chapters 1 and 4]{L87}). The work of
Young \cite{Y01}, \cite{Y02} shows that $d(\lambda)$ gives the
dimension of the irreducible representation of $S_n$ related to
the partition $\lambda\in\Lambda(n)$. A probabilistic proof of
formula (\ref{tableaux}) was given by Greene et al. \cite{GNW79}.
It is based on a construction of a random walk on the cells of the
Young diagram (called also a hook walk). In its first step, a cell
of a Young diagram containing $n$ cells is selected uar. If this
cell has coordinates $(t,s)$, then the walk continues on the cells
$\neq (t,s)$ in the hook of the cell $(t,s)$. A detailed
discussion on the enumerational and algorithmic aspects of formula
(\ref{tableaux}) may be found in \cite[Section 5.1.4]{K98}.

{\it Remark 2.} Our method of proof is based on:
\begin{itemize}
\item{} some combinatorial identities due to Han \cite{H08}(see
Lemmas 1 and 2 in the next section)
 \item{} the saddle point
method in terms of admissibility in the sense of Hayman (see
\cite[Chapter VIII.5]{FS09});
 \item{} the method of moments in probability theory (the
 Fr\'{e}chet-Shohat limit theorem; see, e.g., \cite[Chapter IV, Section 11.4]{L60}).
 \end{itemize}
 Another possible approach to the problem could be built on
 Vershik's ideas for the limit shape of a random Young diagram (see \cite{V96}). In
 this context, one has to deal with a family of probability
 measures $\mu_w, w\in (0,1)$, defined on the set of all
 partitions $\Lambda=\bigcup_{n\ge 0}\Lambda(n)$. The key feature
 in a study of this type is that, under certain conditions, the
 uniform probability measure $\mathbb{P}$ on $\Lambda(n)$ is
 recovered as a conditional distribution $\mu_w(.|\Lambda(n))$.
 Moreover, $\mu_w$ can be constructed as a product measure,
 resulting into mutually independent random part multiplicities $l_j,
 j=1,2,...$ (see their definition given by (\ref{multi})). In the
 context of random partitions, this phenomenon was first
 observed and applied by Fristedt \cite{F93}.  The computation of the expectation of the
 sum of the part sizes $n=n(\lambda)=\lambda_1+\lambda_2+...$ of a partition
 $\lambda=(\lambda_1,\lambda_2,...)\in\Lambda$ with respect to the
 measure $\mu_w$ suggests the proper
 choice of the parameter $w=w(n)$ that approaches $1$ as $n$ becomes large.
  By the above property of the
 conditional distribution $\mu_w(.|\Lambda(n))$, we can consider
 $\lambda\in\Lambda(n)\subset\Lambda$ as selected uar from the set
 $\Lambda(n)$. Dealing with the family of measures $\mu_w$ on the space
 $\Lambda$ avoids
 some technical difficulties because, under this setting, the multiplicities $l_j$ are
 independent. One can now introduce the usual
 scaling for random Young diagrams with $n$
 cells ($A_n=B_n=\sqrt{n}$) in order to obtain their limit shape
 (\ref{limit})-(\ref{shapecur}) in the $(t,s)$-plane as $n\to\infty$. (For more
 details, we refer the reader to \cite{V96} and \cite{B14}).
 Next, we can choose a point at random from the region bounded
 by the lines $t=0, s=0$ and by the curve (\ref{shapecur}). The computation of
 the horizontal and vertical distances between this
 point and the curve is easy. Their sum shows the typical hook length of a
 cell in a random Young diagram for large $n$. Further on, a possible way to identify
 the required limiting distribution is, e.g., to apply some known
 asymptotic results from \cite{P97} or \cite{M151} for the part
 sizes of a random partition $\lambda\in\Lambda(n)$ and of its conjugate partition
 $\lambda^*$. In our further study we prefer, however, to follow an
 approach, which, in our opinion, is more direct: we apply the
 Cauchy formula to a generating function identity and then
 employ the saddle point method for the asymptotic analysis of the
 coefficients representing the moments of the underlying statistic.

\section{Preliminaries}

\subsection{Combinatorial Identities and Generating Functions}

We start with the definition of another partition statistic,
defined on the set $\Lambda(n)$ of all partitions of $n$. We equip
this set with the uniform probability measure and, for any
$\lambda=(\lambda_1,\lambda_2,...,\lambda_k)\in\Lambda(n)$, we
define its $m$-th moment by
\begin{equation}\label{moment}
Y_{m,n}=\sum_{j=1}^k\lambda_j^m,
\end{equation}
In the case $m=0$ the sum in the right-hand side of (\ref{moment})
clearly counts the number of parts in $\lambda$, while, by
(\ref{sigma}), $Y_{1,n}=n$. Further on, $\mathcal{E}(.)$ will
stand for the expectation with respect to the uniform probability
measure on $\Lambda(n)$. From (\ref{moment}) it obviously follows
that
 \begin{equation}\label{expmom}
\mathcal{E}(Y_{m,n})
=\frac{1}{p(n)}\sum_{\lambda\in\Lambda(n)}\sum_j \lambda_j^m,
\quad m=1,2,... .
\end{equation}
Furthermore, the definition of the product probability measure
$\mathbb{P}$ on the space $\Omega(n)$ (see (\ref{pm}) and
(\ref{omega}), respectively) and definition (\ref{z}) of the
random variable $Z_n$ imply that
\begin{equation}\label{expz}
\mathbb{E}(Z_n^m) =\frac{1}{np(n)}\sum_{\lambda\in\Lambda(n)}
\sum_{c\in\lambda} h(\lambda,c), \quad m=1,2,... .
\end{equation}

Our first preliminary facts are related to identities that are due
to Han \cite{H08}. We present the first one in terms of the
expectations introduced by (\ref{expmom}) and (\ref{expz}).

\begin{lemma}
\cite [Corollary 6.5]{H08}. For $m=1,2,...$, we have
$$
\mathbb{E}(Z_n^m)=\frac{1}{n}\mathcal{E}(Y_{m+1,n}).
$$
\end{lemma}

Han \cite[Theorem 6.6]{H08} also obtained an identity for the
generating function of $\mathcal{E}(Y_{m,n})$. It will be the
starting point in our asymptotic analysis. Let
\begin{equation}\label{eul}
g(x)=\prod_{j=1}^\infty (1-x^j)^{-1} =1+\sum_{n=1}^\infty p(n)x^n
\end{equation}
be the Euler partition generating function.

\begin{lemma}
For any $m\ge 1$ and $|x|<1$, we have
$$
1+\sum_{n=1}^\infty p(n)\mathcal{E}(Y_{m,n})x^n =g(x)F_m(x),
$$
where
\begin{equation}\label{efem}
F_m(x) =\sum_{j=1}^\infty\frac{j^m x^j}{1-x^j}.
\end{equation}
\end{lemma}

{\it Remark 3.} In fact, Han \cite{H08} has proved the following
combinatorial identity:
$$
\sum_{\lambda\in\Lambda(n)}\sum_{c\in\lambda} h^m(\lambda,c)
=\sum_{\lambda\in\Lambda(n)}\sum_j \lambda_j^{m+1}.
$$
Lemma 1 translates it in terms of the expectations (\ref{expmom})
and (\ref{expz}).

\subsection{Analytic Combinatorics Background: Meinardus Theorem
on Weighted Partitions and Hayman Admissibility of Generating
Functions}

Lemma 2 implies that the coefficient $p(n)\mathcal{E}(Y_{m,n})$ of
$x^n$ can be expressed by a Cauchy integral whose integrand
contains the product $g(x)F_m(x)$. Its behavior heavily depends on
the properties of the partition generating function $g(x)$ whose
infinite product representation (\ref{eul}) shows that its main
singularity is at the point $x=1$ (see also \cite[Chapter
5]{A76}). A complete asymptotic expansion of
$\mathcal{E}(Y_{m,n})$, for fixed and odd $m>1$, was obtained by
Grabner et al. \cite{GKW14}. We need here only the main asymptotic
term of $\mathcal{E}(Y_{m,n})$, $m=1,2,...$, which can be obtained
using the Hayman-Meinardus version of the saddle point method.
Here, we will give a brief account on this subject. This approach
was also demonstrated in \cite[Sections 5 and 6]{M151}.

We start with some general remarks on the Hardy-Ramanujan formula
(\ref{hr}). The Cauchy integral stemming from (\ref{eul}) and
representing $p(n)$ can be analyzed with the aid of the Hardy,
Ramanujan and Rademacher's circle method developed in \cite{HR18}
and \cite{R37}. Subsequent generalizations of these results are
related to extensions of the class of generating functions, for
which an infinite product representation similar to (\ref{eul}) is
valid. An important result in this direction is due to Meinardus
\cite{M54} (see also \cite[Chapter 6]{A76}) who obtained the
asymptotic of the Taylor coefficients of infinite products of the
form:
\begin{equation}\label{ip}
\prod_{j=1}^\infty (1-x^j)^{-b_j}
\end{equation}
under certain general assumptions on the sequence of non-negative
numbers $\{b_j\}_{j\ge 1}$. Meinardus' approach is based on
considering the Dirichlet generating series
\begin{equation}\label{diri}
D(z)=\sum_{j=1}^\infty b_j j^{-z}, \quad z=u+iv, \quad
u,v\in\mathbb{R}.
\end{equation}
We briefly describe here Meinardus' assumptions avoiding their
precise statements as well as some extra notations and concepts.

\begin{itemize}
\item{$(M_1)$} The first assumption ($M_1$) specifies the domain
$\mathcal{H}=\{z=u+iv:u\ge -C_0\}, 0<C_0<1,$ in the complex plane,
in which $D(z)$ has an analytic continuation.

\item{$(M_2)$} The second one ($M_2$) is related to the asymptotic
behavior of $D(z)$, whenever $|v|\to\infty$. A function of the
complex variable $z$ which is bounded by $O(|\Im{(z)}|^{C_1}),
0<C_1<\infty$, in certain domain of the complex plane is called a
function of finite order. Meinardus' second condition requires
that $D(z)$ is of finite order in the whole domain $\mathcal{H}$.

\item{$(M_3)$} Finally, the Meinardus third condition ($M_3$)
implies a bound on the ordinary generating function of the
sequence $\{b_j\}_{j\ge 1}$. It can be stated in a way simpler
than the Meinardus' original expression by the inequality:
$$
\sum_{j\ge1} b_j e^{-j\alpha}\sin^2{(\pi jy)}\ge
C_2\alpha^{-\epsilon_1}, \quad
0<\frac{\alpha}{2\pi}<|y|<\frac{1}{2},
$$
for sufficiently small $\alpha$ and some constants $C_2,
\epsilon_1>0$ ($C_2=C_2(\epsilon_1)$) (see \cite{GSE08}).
\end{itemize}

It is known that the Euler partition generating function $g(x)$
(which is obviously of the form (\ref{ip}) satisfies the Meinardus
scheme of conditions ($M_1$) - ($M_3$) (see, e.g., \cite[Theorem
6.3]{A76}).

In the asymptotic analysis of the Cauchy integral stemming from
Lemma 2 we shall use a variant of the saddle point method given by
Hayman's theorem for admissible power series \cite{H56} (for more
details, see, e.g., \cite[Chapter VIII.5]{FS09}). To present
Hayman's idea and show how it can be applied to the proof of
Theorem 1, we need to introduce some auxiliary notations.

We consider here a function $G(x)=\sum_{n=1}^\infty G_n x^n$ that
is analytic for $|x|<\rho$, $0<\rho<\infty$. For $0<r<\rho$, we
let
\begin{equation}\label{a}
a(r)=r\frac{G^\prime(r)}{G(r)},
\end{equation}

\begin{equation}\label{b}
b(r)=r\frac{G^\prime(r)}{G(r)}
+r^2\frac{G^{\prime\prime}(r)}{G(r)}
-r^2\left(\frac{G^\prime(r)}{G(r)}\right)^2.
\end{equation}

In the statement of the Hayman's result, we shall use the
terminology given in \cite[Chapter VIII.5]{FS09}. We assume that
$G(x)>0$ for $x\in (R_0,\rho)\subset (0,\rho)$ and satisfies the
following three conditions.

\begin{itemize}
\item {\it Capture condition.} $\lim_{r\to\rho} a(r)=\infty$ and
$\lim_{r\to\rho} b(r)=\infty$.

\item {\it Locality condition.} For some function
$\delta=\delta(r)$ defined over $(R_0,\rho)$ and satisfying
$0<\delta<\pi$, one has
$$
G(r e^{i\theta})\sim G(r) e^{i\theta a(r)-\theta^2 b(r)/2}
$$
as $r\to\rho$, uniformly for $|\theta|\le\delta(r)$.

\item {\it Decay condition.}
$$
G(r e^{i\theta})= o\left(\frac{G(r)}{\sqrt{b(r)}}\right)
$$
as $r\to\rho$, uniformly for $\delta(r)\le |\theta|<\pi$.
\end{itemize}

\begin{definition}
A function $G(x)$ satisfying the capture, locality and decay
conditions is called a Hayman admissible function.
\end{definition}

{\bf Hayman Theorem.} Let $G(x)$ be a Hayman admissible function
and $r=r_n$ be the unique solution in the interval $(R_0,\rho)$ of
the equation
\begin{equation}\label{sp}
a(r)=n.
\end{equation}
Then the Taylor coefficients of $G(x)$ satisfy, as $n\to\infty$,
$$
G_n\sim\frac{G(r_n)}{r_n^n\sqrt{2\pi b(r_n)}}
$$
with $a(r_n)$ and $b(r_n)$ given by (\ref{a}) and (\ref{b}),
respectively.

\section{Proof of Theorem 1}

The proof is divided into two parts.

(A) Proof of Hayman admissibility for $g(x)$.

(B) Obtaining an asymptotic estimate for the Cauchy integral
stemming from Lemma 2.

\subsection{Part (A)}

This part of the proof follows the same line of reasoning given in
\cite[pp. 338-340]{M151}. For the sake of completeness, we present
it here.

First, we need to show how Hayman theorem can be applied to find
the asymptotic behavior of the Taylor coefficients of the
partition generating function $g(x)$. Since in (\ref{ip}) we have
$b_j=1, j\ge 1$, the Dirichlet generating series (\ref{diri}) is
$D(z)=\zeta(z)$, where $\zeta$ denotes the Riemann zeta function:
$\zeta(z)=\sum_{j=1}^\infty j^{-z}, z=u+iv$. We set in (\ref{sp})
$r=r_n=e^{-d_n}, d_n>0,$ where $d_n$ is the unique solution of the
equation
\begin{equation}\label{expsp}
a(e^{-d_n})=n.
\end{equation}
Granovsky et al. \cite{GSE08} showed that the first two Meinardus
conditions imply that the unique solution of (\ref{expsp}) has the
following asymptotic expansion:
\begin{eqnarray}\label{d}
& & d_n=\sqrt{\zeta(2)/n}+\frac{\zeta(0)}{2n} +O(n^{-1-\beta})
\nonumber \\
& & =\frac{\pi}{\sqrt{6n}} -\frac{1}{4n} +O(n^{-1-\beta}),
\end{eqnarray}
where $\beta>0$ is fixed constant (here we have also used that
$\zeta(0)=-1/2$; see \cite[Chapter 23.2]{AS65}. We also notice
that (\ref{b}) and (\ref{d}) imply that
\begin{equation}\label{asb}
b(e^{-d_n}) =2\zeta(2)d_n^{-3} +O(d_n^{-2})
\sim\frac{2\sqrt{6}}{\pi}n^{3/2}.
\end{equation}
(This is a particular case of Lemma 2.2 from \cite{M13} with
$D(z)=\zeta(z)$.) Hence, by (\ref{expsp}) and (\ref{asb}),
$a(e^{-d_n})\to\infty$ and $b(e^{-d_n})\to\infty$ as $n\to\infty$,
that is, Hayman's "capture" condition is satisfied with
$r=r_n=e^{-d_n}$. To show next that Hayman's "decay" condition is
satisfied, we set
\begin{equation}\label{del}
\delta_n =\frac{d_n^{4/3}}{\omega(n)}
=\frac{\pi^{4/3}}{(6n)^{2/3}\omega(n)}
\left(1+O\left(\frac{1}{\sqrt{n}}\right)\right),
\end{equation}
where $\omega(n)$ is a function satisfying $\omega(n)\to\infty$ as
$n\to\infty$ arbitrarily slowly and in the second equality we have
used (\ref{d}). We can apply now an estimate for
$|g(e^{-d_n+i\theta})|$ established in a general form in
\cite[Lemma 2.4]{M13} with the aid of all three Meinardus
conditions. In our particular case it states that there exist two
positive constants $c_0$ and $\epsilon_0$, such that, for
sufficiently large $n$,
\begin{equation}\label{err}
|g(e^{-d_n+i\theta})|\le g(e^{-d_n}) e^{-c_0 d_n^{-\epsilon_0}}
\end{equation}
uniformly for $\delta_n\le |\theta|<\pi$. This, in combination
with (\ref{d}) and (\ref{asb}) implies that
$|g(e^{-d_n+i\theta})|= o(g(e^{-d_n})/\sqrt{b(e^{-d_n})})$
uniformly in the same range of $\theta$, which is just Hayman's
"decay" condition. Finally, by Lemma 2.3 of \cite{M13},
established with the aid of Meinardus' conditions ($M_1$) and
($M_2$), Hayman's "locality" condition is also satisfied by
$g(x)$. In fact, this lemma implies in the particular case
$D(z)=\zeta(z)$ that
\begin{equation}\label{main}
e^{-i\theta n}\frac{g(e^{-d_n+i\theta})}{g(e^{-d_n})}
=e^{-\theta^2
b(e^{-d_n})/2}\left(1+O\left(\frac{1}{\omega^3(n)}\right)\right)
\end{equation}
uniformly for $|\theta|\le\delta_n$, where $b(e^{-d_n})$ and
$\delta_n$ are determined by (\ref{asb}) and (\ref{del}),
respectively. Hence, all conditions of Hayman's theorem hold, and
we can apply it with $G_n=p(n), G(x)=g(x), r_n=e^{-d_n}$ and
$\rho=1$ to find that
\begin{equation}\label{pn}
p(n)\sim \frac{e^{nd_n} g(e^{-d_n})}{\sqrt{2\pi b(e^{-d_n})}},
\quad n\to\infty.
\end{equation}

{\it Remark 4.} To show that formula (\ref{pn}) yields (\ref{hr}),
one has to replace (\ref{d}) and (\ref{asb}) in the right-hand
side of (\ref{pn}). The asymptotic of $g(e^{-d_n})$ is determined
by a general lemma due to Meinardus \cite{M54} (see also
\cite[Lemma 6.1]{A76}). Since $\zeta(0)=-1/2$ and
$\zeta^\prime(0)=-\frac{1}{2}\log{(2\pi)}$ (see \cite[Chapter
23.2]{AS65}), in the particular case of $g(e^{-d_n})$, this lemma
implies that
\begin{eqnarray}
& & g(e^{-d_n})
=\exp{(\zeta(2)d_n^{-1}-\zeta(0)\log{d_n}+\zeta^\prime(0)
+O(d_n^{c_1}))} \nonumber \\
& & =\exp{\left(\frac{\pi^2}{6d_n}+\frac{1}{2}\log{(2\pi)}
+O(d_n^{c_1})\right)}, \quad n\to\infty , \nonumber
\end{eqnarray}
where $0<c_1<1$. The rest of the computation leading to (\ref{hr})
is based on simple algebraic manipulations and cancellations.

\subsection{Part (B)}

We are now ready to apply the Cauchy coefficient formula to the
generating function identity of Lemma 2. We use the circle
$x=e^{-d_n+i\theta}, -\pi<\theta\le\pi$, as a contour of
integration and obtain, for any fixed $m\ge 1$, that
$$
p(n)\mathcal{E}(Y_{m,n}) =\frac{e^{n d_n}}{2\pi}\int_{-\pi}^\pi
g(e^{-d_n+i\theta}) F_m(e^{-d_n+i\theta}) e^{-i\theta n} d\theta.
$$
Then, we break up the range of integration as follows:
\begin{equation}\label{onetwo}
p(n)\mathcal{E}(Y_{m,n}) =J_1(m,n)+J_2(m,n),
\end{equation}
where
\begin{equation}\label{one}
J_1(m,n)=\frac{e^{n d_n}}{2\pi}\int_{-\delta_n}^{\delta_n}
g(e^{-d_n+i\theta}) F_m(e^{-d_n+i\theta}) e^{-i\theta n} d\theta,
\end{equation}
\begin{equation}\label{two}
J_2(m,n)=\frac{e^{n d_n}}{2\pi}\int_{\delta_n<|\theta|\le\pi}
g(e^{-d_n+i\theta}) F_m(e^{-d_n+i\theta}) e^{-i\theta n} d\theta,
\end{equation}
where $m=1,2,...$, and $\delta_n$ is defined by (\ref{del}).

To estimate $J_2(m,n)$, we notice that, for fixed $m$, by the
definition of Riemann integrals, (\ref{efem}) and (\ref{d}),
\begin{eqnarray}\label{ineq}
& & |F_m(e^{-d_n+i\theta})| =|\sum_{j=1}^\infty\frac{j^m e^{-j
d_n+ij\theta}}{1-e^{-j d_n+ij\theta}}| \nonumber \\
& & \le\sum_{j=1}^\infty\frac{j^m e^{-j d_n}}{1-e^{-j d_n}}
=\sum_{j\ge 1}\frac{(j d_n)^m e^{-j d_n}}{1-e^{-j d_n}} d_n^{-m-1}
\\
& & \sim d_n^{-m-1} \int_0^\infty\frac{u^m e^{-u}}{1-e^{-u}} du
=O(d_n^{-m-1})= O(n^{(m+1)/2}). \nonumber
\end{eqnarray}
The last two equalities follow from the estimate
\begin{equation}\label{dminus}
d_n^{-1}=\frac{\sqrt{6n}}{\pi} +O(1),
\end{equation}
which is a simple consequence of (\ref{d}). Combining (\ref{two})
- (\ref{dminus}), (\ref{err}) and (\ref{asb}) with the asymptotic
equivalence (\ref{pn}) for the numbers $p(n)$, we obtain
\begin{eqnarray}\label{twolast}
& & |J_2(m,n)|\le\frac{e^{n d_n}}{2\pi}\int_{\delta_n\le
|\theta|<\pi} |g(e^{-d_n+i\theta})| |F_m(e^{-d_n+i\theta})|
d\theta
\nonumber \\
& & =O(e^{n d_n} g(e^{-d_n}) n^{(m+1)/2} e^{-c_0
d_n^{-\epsilon_0}}) \nonumber \\
& & O\left(\frac{e^{n d_n} g(e^{-d_n})}{\sqrt{b(e^{-d_n})}}
n^{(m+1)/2} n^{3/4} e^{-c_2 n^{\epsilon_0/2}}\right) \nonumber \\
& & =O(p(n) n^{(m+1)/2} n^{3/4} e^{-c_2 n^{\epsilon_0/2}}) =o(p(n)
n^{(m+1)/2}),
\end{eqnarray}
where $c_2>0$.

For the asymptotic estimate of $J_1(m,n)$, we need to expand
$F_m(x)$ around the point $x=e^{-d_n}$. Thus, for any fixed
$m=1,2,...$ and uniformly for any $|\theta|\le\delta_n$, we have
\begin{eqnarray}\label{firstef}
& & F_m(e^{-d_n+i\theta}) =F_m(e^{-d_n})
+O\left(|\theta|\frac{d}{dx}F_m(x)|_{x=e^{-d_n}}\right) \nonumber \\
& & =F_m(e^{-d_n})
+O\left(\delta_n\frac{d}{dx}F_m(x)|_{x=e^{-d_n}}\right).
\end{eqnarray}
As previously, we can consider the sum representing
$F_m(e^{-d_n})$ as a Riemann sum. So, for large $n$, we can
replace it by the value of the corresponding integral (see, e.g.,
\cite[Section 27.1]{AS65}). Hence, by (\ref{efem}), (\ref{d}) and
(\ref{dminus}), we have
\begin{eqnarray}\label{secef}
& & F_m(e^{-d_n}) =\sum_{j=1}^\infty\frac{j^m e^{-j d_n}}{1-e^{-j
d_n}} \sim d_n^{-m-1} \int_0^\infty\frac{u^m e^{-u}}{1-e^{-u}} du
\nonumber \\
& & =d_n^{-m-1}
m!\zeta(m+1)\sim\left(\frac{n}{\zeta(2)}\right)^{(m+1)/2}
m!\zeta(m+1), \quad m=1,2,... .
\end{eqnarray}
In the same way, we can estimate the first derivative of $F_m$:
\begin{eqnarray}\label{deref}
& & \frac{d}{dx} F_m(x)|_{x=e^{-d_n}}
=\sum_{j=1}^\infty\frac{j^{m+1} e^{-j d_n} e^{d_n}}
{(1-e^{-jd_n})^2} \nonumber \\
& & \sim d_n^{-m-2}\int_0^\infty\frac{u^{m+1}e^{-u}}
{(1-e^{-u})^2} du =O(d_n^{-m-2}) =O(n^{(m+2)/2}),
\end{eqnarray}
where the last $O$-estimate follows from (\ref{dminus}). Hence, by
(\ref{del}) and (\ref{deref}), the error term in (\ref{firstef})
becomes
\begin{eqnarray}\label{errt}
& & O(\delta_n n^{(m+2)/2}) =O(n^{(m+1)/2} n^{1/2}
n^{-2/3}/\omega(n)) \nonumber \\
& & =O(n^{(m+1)/2} n^{-1/6}/\omega(n)) =o(n^{(m+1)/2}).
\end{eqnarray}
Consequently, (\ref{firstef}) - (\ref{errt}) imply that
$$
F_m(e^{-d_n+i\theta}) =\left(\frac{n}{\zeta(2)}\right)^{(m+1)/2}
m!\zeta(m+1)+o(n^{(m+1)/2}).
$$
Inserting this estimate and (\ref{main}) into (\ref{one}) and
applying the asymptotic of the partition function $p(n)$ from
(\ref{pn}), we obtain
\begin{eqnarray}\label{oneest}
& & J_1(m,n) =\frac{e^{n d_n}g(e^{-d_n})}{2\pi}
\int_{-\delta_n}^{\delta_n}
\frac{g(e^{-d_n+i\theta})}{g(e^{-d_n})} \nonumber \\
& & \times\left(\left(\frac{n}{\zeta(2)}\right)^{(m+1)/2}
m!\zeta(m+1) +o(n^{(m+1)/2})\right) e^{-i\theta n} d\theta
\nonumber \\
& & =\frac{e^{n d_n}g(e^{-d_n})}{2\pi}
\left(\left(\frac{n}{\zeta(2)}\right)^{(m+1)/2} m!\zeta(m+1)
+o(n^{(m+1)/2})\right) \nonumber \\
& & \times\int_{-\delta_n}^{\delta_n} e^{-\theta^2 b(e^{-d_n})/2}
(1+1/\omega^3(n))d\theta \nonumber \\
& & =\frac{e^{n d_n}g(e^{-d_n})}{2\pi\sqrt{b(e^{-d_n})}}
\left(\left(\frac{n}{\zeta(2)}\right)^{(m+1)/2} m!\zeta(m+1)
+o(n^{(m+1)/2})\right) \nonumber \\
& & \times
\int_{-\delta_n\sqrt{b(e^{-d_n})}}^{\delta_n\sqrt{b(e^{-d_n})}}
e^{-y^2/2}dy \nonumber \\
& & \sim \frac{e^{n d_n}g(e^{-d_n})}{2\pi\sqrt{b(e^{-d_n})}}
\left(\frac{n}{\zeta(2)}\right)^{(m+1)/2} m!\zeta(m+1)
\int_{-\infty}^\infty e^{-y^2/2}dy \nonumber \\
& & \sim \frac{e^{n d_n}g(e^{-d_n})}{\sqrt{2\pi b(e^{-d_n})}}
\left(\frac{n}{\zeta(2)}\right)^{(m+1)/2} m!\zeta(m+1) \nonumber
\\
& & \sim p(n) \left(\frac{n}{\zeta(2)}\right)^{(m+1)/2}
m!\zeta(m+1).
\end{eqnarray}
In the first asymptotic equivalence we have used (\ref{asb}) and
(\ref{del}) in order to get
$$
\delta_n\sqrt{b(e^{-d_n})}
\sim\frac{\pi^{5/6}\sqrt{2}}{6^{1/6}\omega(n)} n^{1/12}\to\infty
$$
if $\omega(n)\to\infty$ as $n\to\infty$ not too fast, so that
$n^{1/12}/\omega(n)\to\infty$. Combining (\ref{onetwo}) -
(\ref{two}), (\ref{twolast}) and (\ref{oneest}), we get
$$
p(n)\mathcal{E}(Y_{m,n})
=p(n)\left(\frac{n}{\zeta(2)}\right)^{(m+1)/2} m!\zeta(m+1)
+o(n^{(m+1)/2}p(n)),
$$
which in turn shows that
\begin{equation}\label{mom}
\mathcal{E}(Y_{m,n}) \sim\left(\frac{n}{\zeta(2)}\right)^{(m+1)/2}
m!\zeta(m+1), \quad m=1,2,... .
\end{equation}

Now, we recall Lemma 1, (\ref{expmom}) and (\ref{expz}). Combining
these observations with (\ref{mom}), we obtain
\begin{eqnarray}
& & \mathbb{E}(Z_n^m) =\frac{1}{n p(n)}\sum_{\lambda\in\Lambda(n)}
\sum_j \lambda_j^{m+1} =\frac{1}{n}\mathcal{E}(Y_{n,m+1})
\nonumber
\\
& & \sim \frac{1}{n}\left(\frac{n}{\zeta(2)}\right)^{m/2+1}
(m+1)!\zeta(m+2), \quad m=1,2,... . \nonumber
\end{eqnarray}
Consequently,
\begin{eqnarray}
& & \lim_{n\to\infty}
\mathbb{E}\left(\left(\sqrt{\frac{\zeta(2)}{n}} Z_n\right)^m
\right) =\lim_{n\to\infty}
\mathbb{E}\left(\left(\frac{\pi}{\sqrt{6n}}
Z_n\right)^m \right) \nonumber \\
& & =\frac{(m+1)!\zeta(m+2)}{\zeta(2)}
=\frac{6}{\pi^2}\int_0^\infty\frac{u^{m+1}}{e^u-1} du, \quad
m=1,2,... . \nonumber
\end{eqnarray}
Hence, the Frechet-Shohat limit theorem \cite[Chapter IV, Section
11.4]{L60} implies that the sequence $\{\pi Z_n/\sqrt{6n}\}_{n\ge
1}$ converges in distribution to a random variable with
probability density function $\frac{6u}{\pi^2(e^u-1)}, u\ge 0$,
and zero elsewhere, which completes the proof of the theorem.


\begin{thebibliography}{15}

\bibitem{AS65}
M. Abramovitz and I. A. Stegun, {\it Handbook of Mathemathical
Functions with Formulas, Graphs and Mathematical Tables}, Dover
Publ. Inc. (New York, 1965).

\bibitem{A76}
G. E. Andrews, {\it The Theory of Partitions, Encyclopedia Math.
Appl. 2}, Addison-Wesley (Reading, MA, 1976).

\bibitem{B14}
L. V. Bogachev, Unified derivation of the limit shape
formultiplicative ensembles of random integer partitions with
equiweighted parts, {\em Random Struct. Alg.} {\bf 47} (2015),
227-266.

\bibitem{CPSW99}
S. Corteel, B. Pittel, C. D. Savage, H. S. Wilf, On the
multiplicity of parts in a random partition, {\em Random Stuct.
Alg.} {\bf 14} (1999), 185-197.

\bibitem{EL41}
P. Erd\"{o}s and J. Lehner, The distribution of the number of
summands in the partitions of a positive integer, {\em Duke Math.
J.} {\bf 8} (1941), 335-345.

\bibitem{ES84}
P. Erd\"{o}s  and M. Szalay, On the statistical theory of
partitions, {\em In: Topics in Classical Number Theory, vol. I (G.
Halasz ed.)}, North-Holland, Amsterdam, pp. 397-450, 1984.

\bibitem{FS09}
P. Flajolet and R. Sedgewick, {\it Analytic Combinatorics},
Cambridge Univ. Press (Cambridge, 2009).

\bibitem{FRT54}
J. S. Frame, G. de B. Robinson and R.M. Thrall, The hook graphs of
the symmetric group, {\em Canad. J. Math.} {\bf 6} (1954),
316-324.

\bibitem{F93}
B. Fristedt, The structure of random partitions of large integers,
{\em Trans. Amer. Math. Soc.}  {\bf 337} (1993), 703-735.

\bibitem{GKW14}
P. Grabner, A. Knopfmacher and S. Wagner, A general asymptotic
scheme for the analysis of partition statistics, {\em Combin.
Probab. Comput.} {\bf 23} (2014), 1057-1086.

\bibitem{GSE08}
B. L. Granovsky, D. Stark and M. Erlihson, Meinardus' theorem on
weighted partitions: Extensions and a probailistic proof. {\em
Adv. Appl. Math.} {\bf 41} (2008) 307-328.

\bibitem{GNW79}
C. Greene, A. Nijenhuis and H. Wilf, A probabilistic proof of a
formula for the number of Young tableaux of a given shape, {\em
Adv. Math.} {\bf 31} (1979), 104-109.

\bibitem{H08}
G.-N. Han, An explicit expansion formula for the powers of the
Euler product in terms of partition hook lengths, arXiv: 0804-1849
v3 [math.CO] (2008).

\bibitem{HR18}
G. H. Hardy and S. Ramanujan, Asymptotic formulae in combinatory
analysis, {\em Proc. London Math. Soc.} {\bf 17(2)} (1918),
75-115.

\bibitem{H56}
W. K. Hayman, A generalization of Stirling's formula, {\em J.
Reine Angew. Math.} {\bf 196} (1956), 67-95.

\bibitem{K98}
D. E. Knuth, {\it The Art of Computer Programming, Vol. 3: Sorting
and Searching, 2nd ed.}, Addison-Wesley-Longman (Reading, MA,
1998).

\bibitem{L87}
W. Ledermann, {\it Introduction to Group Chatracters, 2nd ed.},
Cambridge Univ. Press (Cambridge, 1987).

\bibitem{L60}
M. Lo\'{e}ve, {\it Probability Theory, 2nd ed.}, Van Nostrand,
(Princeton, NJ, 1960).

\bibitem{M54}
G. Meinardus, Asymptotische Aussagen \"uber Partitionen. {\em
Math. Z.} {\bf 59} (1954), 388-398.

\bibitem{M05}
L. Mutafchiev, On the maximal multiplicity of parts in a random
integer partition, {\em Ramanujan J.} {\bf 9} (2005), 305-316.

\bibitem{M13}
L. Mutafchiev, The size of the largest part of random weighted
partitions of large integers, {\em Combin. Probab. Comput.} {\bf
22} (2013), 433-454.

\bibitem{M151}
L. Mutafchiev, Sampling part sizes of random integer partitions,
{\em Ramanujan J.} {\bf 37} (2015), 329-343.

\bibitem{M152}
L. Mutafchiev, Sampling parts of random integer partitions: a
probabilistic and asymptotic analysis, {\em Pure Math. Appl.} {\bf
25} (2015), 79-95.

\bibitem{P97}
B. Pittel, On a likely shape of the random Ferrers diagram, {\em
Adv. Appl. Math.} {\bf 18} (1997), 432-488.

\bibitem{R37}
H. Rademacher, On the partition function $p(n)$, {\em Proc. London
Math. Soc.} {\bf 43} (1937), 241-254.

\bibitem{ST771}
M. Szalay and P. Tur\'{a}n, On some problems of the statistical
theory of partitions with applications to characters of the
symmetric group, I, {\em Acta Math. Acad. Sci. Hungar.} {\bf 29}
(1977), 361-379.

\bibitem{ST772}
M. Szalay and P. Tur\'{a}n, On some problems of the statistical
theory of partitions with applications to characters of the
symmetric group, II, {\em Acta Math. Acad. Sci. Hungar.} {\bf 29}
(1977), 381-392.

\bibitem{ST78}
M. Szalay and P. Tur\'{a}n, On some problems of the statistical
theory of partitions with applications to characters of the
symmetric group, III, {\em Acta Math. Acad. Sci. Hungar.} {\bf 32}
(1978), 129-155.

\bibitem{T52}
H.N.V. Temperley, Statistical mechanics and the partition of
numbers II. The form of crystal surfaces, {\em Math. Proc.
Cambridge. Philos. Soc.} {\bf 48} (1952), 683-697.

\bibitem{V96}
A.M. Vershik, Statistical mechanics of combinatorial partitions
and their limit shapes, {\em Funct. Anal. Appl.} {\bf 30} (1996),
90-105.

\bibitem{W11}
S. Wagner, Limit distributions of smallest gap and largest
repeated part in integer partitions, {\em Ramanujan J.} {\bf 25}
(2011), 229-246.

\bibitem{Y01}
A. Young, On quantitative substitutional analysis, I, {\em Proc.
London Math. Soc.} {\bf 1(32)} (1901), 384-404.

\bibitem{Y02}
A. Young, On quantitative substitutional analysis, II, {\em Proc.
London Math. Soc.} {\bf 1(34)} (1902), 202-208.

\end{thebibliography}
\end{document}